\documentclass[11pt]{amsart}
\usepackage{graphicx}

\newtheorem{thm}{Theorem}[section]

\theoremstyle{definition}

\theoremstyle{remark}
\newtheorem{rem}[thm]{Remark}
\numberwithin{equation}{section}

\begin{document}

\begin{center}
{\Large Algebraic Sum of Unbounded Normal Operators and the Square
Root Problem of Kato}
\end{center}

\vspace{0.9cm}

\begin{center}
{\bf By Toka Diagana}\footnotetext[1] { Howard University,
Department of Mathematics, 2441 6th Street, N.W - Washington, D.C
20059 - USA / E-Mail: tdiagana@howard.edu}
\end{center}

\vspace{0.7cm}

\noindent {\bf Abstract}:{\small We prove that the algebraic sum
of unbounded normal operators satisfies the square root problem of
Kato under appropriate hypotheses. As application, we consider
perturbed Schr\"{o}dinger operators}.


\noindent 

\vspace{0.2cm}

\noindent 

\section{Introduction}
Let $C$ be a normal operator (not necessarily bounded) in a
(complex) Hilbert space ${\mathcal H}$. Using the spectral theorem
for unbounded normal operators, it is well-known that $C$ 
can be expressed as
\begin{equation}
C = C_{1} - i C_{2} 
\end{equation}
with $C_{1}$, $C_{2}$ 
unbounded selfadjoint operators on ${\mathcal H}$ (see,  e.g [ 13 , pp. 348-355]). If one supposes $C_{1}$, $C_{2}$ to be
nonnegative operators, then $\; i C \;$ is m-accretive (see,  e.g, [ 12, Corollary 4.4, p. 15]).

Let $A$, $B$ be normal operators on
${\mathcal H}$ . Recall the algebraic sum 
$S = A+B$ of $A$ and $B$ is defined as
\begin{equation}
\forall u \in D(S) = D(A) \cap D(B), \; \; \; S u = A u + B u
\end{equation}
In this paper we are concerned with the square root problem of Kato
for the sum of operators $A$ and $B$.

In section 2, we prove that the algebraic sum $S$ 
satisfies the square root problem of Kato under suitable
hypotheses, that is:
\begin{equation}
D(S^{\frac{1}{2}}) = D(A^{\frac{1}{2}}) \cap D( B^{\frac{1}{2}}) =
D(S^{* {\frac{1}{2}}})
\end{equation}
Also, since the algebraic sum $S$ is is not always defined (see, e.g., [ 5, 8 ]).
We shall define a "generalized" sum $A \oplus B$ of $A$ and $B$. One can
then prove that such a "generalized" sum satisfies the
square root problem of Kato under suitable hypotheses. As application we
shall consider perturbed Schr\"{o}dinger operators.

Recall that more details about the well-known square root problem of Kato
can be found in [1, 4, 7, 9, 11].

Throughout this paper we assume that $A$ and $B$ can be decomposed
as $A = A_{1} - i A_{2}$ and $B = B_{1} - i B_{2}$. We denote by $\Omega = \Omega(A) \cap \Omega(B)$ where
$$\Omega(A) = D( |A|^{\frac{1}{2}}) = D(A_{1}^{\frac{1}{2}}) \cap
D(A_{2}^{\frac{1}{2}}), \; \; \; \Omega(B) = D( |B|^{\frac{1}{2}})
= D(B_{1}^{\frac{1}{2}}) \cap D(B_{2}^{\frac{1}{2}})$$ 

The following assumptions will be made

\vspace{0.3cm}

\noindent $ (H_{1}) \; \; \; \; A_{k}, \; \; B_{k} \; \; \; \mbox{are nonnegative}\;
(k = 1, 2)$

\vspace{0.2cm}

\noindent $(H_{2}) \; \; \; \; \overline{D(A) \cap D(B)} = {\mathcal H}
$
\vspace{0.2cm}

\noindent $(H_{3}) \; \; \; \; \mbox{there exists} \; \; c > 0, \;
\; \; < A_{2} u , u
> \; \leq \; c \; < A_{1} u , u >, \; \; \forall u \in \Omega(A)$

\vspace{0.2cm}

\noindent $(H_{4}) \; \; \; \; \mbox{there exists} \; \; c' > 0,
\; \; \; < B_{2} u , u
> \; \leq \; c' \; < B_{1} u , u >, \; \; \forall u \in
\Omega(B)$

\vspace{0.2cm}

\noindent $(H_{5}) \; \; \; \; \overline{\Omega} = {\mathcal H}
$
\vspace{0.2cm}

\noindent $ (H_{6}) \; \;\; \; \Omega$ is closed in the interpolation space
$[ \Omega \; , \; {\mathcal H}]_{\frac{1}{2}}$

\vspace{0.3cm}

\noindent Consider the sesquilinear forms associated with $A$
and $B$:
$$\phi(u,v) \; = \; < \; A_{1}^{\frac{1}{2}}u ,
A_{1}^{\frac{1}{2}} v \; > - i < \; A_{2}^{\frac{1}{2}} u ,
A_{2}^{\frac{1}{2}}v \; >, \; \; \; \forall \; u,v \in \Omega(A)$$

$$\psi(u,v) \; = \; < \; B_{1}^{\frac{1}{2}}u ,
B_{1}^{\frac{1}{2}} v \; > - i < \; B_{2}^{\frac{1}{2}} u ,
B_{2}^{\frac{1}{2}}v \; > \; \; \; \forall \; u,v \in \Omega(B)$$

Set
\begin{equation}
\xi (u,v) = \phi(u,v) +
\psi(u,v), \; \; \; \; \; \forall u,v \in \Omega
\end{equation}
According to Bivar-Weinholtz and Lapidus (see,  e.g., [2 , pp. 451])
the "Generalized" sum $A \oplus B$ of $A$ and $B$ is defined with
the help of the sesquilinear form $\xi$ as follows: $u \in D(A \oplus B)$ iff $\; v
\longrightarrow \xi(u,v)$ is continuous for the
${\mathcal{H}}$-Topology, and $(A \oplus B)u$ defined to be the vector of
${\mathcal{H}}$ given by the Riesz-Representation theorem
\begin{equation}
< \; (A \oplus B) u , v \; > \; = \; \xi (u,v) \; \; \forall v \in
\Omega
\end{equation}
Applying $(H_{2})$ to (1.4), we see that $\xi$
admits the following representation
\begin{equation}
\xi(u,v) \; = \; < \; (A+B) u , v \; >, \; \; \forall u \in D(A)
\cap D(B), \; \; \forall v \in \Omega
\end{equation}

\section{Square Root Problem of Kato}
In this section we show the algebraic
sum $S = A+B$ satisfies the square root problem of Kato under suitable conditions.
We also show the same conclusion still holds if we consider the square root problem of Kato for the "generalized" sum defined above.

We have

{\thm Let $A = A_{1} - i A_{2}$ and $B = B_{1} - i B_{2}$ be unbounded normal operators on 
${\mathcal{H}}$. Assume that assumptions
$(H_{1}), \; (H_{2}), \; (H_{3}), \; \;  \mbox{and} \;
\; (H_{4})$ are satisfied and that the operator $\overline{A+B}$ is 
maximal. Then we have
$$D((\overline{A+B})^{\frac{1}{2}}) = D(A^{\frac{1}{2}}) \cap D(B^{\frac{1}{2}}) =D((\overline{A+B})^{*
\frac{1}{2}})$$}

\noindent Proof. Consider the sesquilinear form $\xi = \phi +
\psi$ given by (1.6). Let $\Omega_{\xi} = ( \Omega , \| \; . \;
\|_{\xi})$ be the pre-Hilbert space $\Omega$ equipped with the
inner product given as $$< \; u,v \; >_{\xi} \; = \; < \; u,v \;
>_{{\mathcal{H}}}\; + \; \Re e \xi(u,v), \; \; \; \; \forall \; u,v \in \Omega$$
Since the sum form  $A_{1} \oplus B_{1}$ is a nonnegative
selfadjoint operator. It easily follows that $\Omega_{\xi}$ is a
Hilbert space, and therefore $\xi$ is a closed sesquilinear form.
Moreover, $D(\xi) = \Omega$ is dense on ${\mathcal{H}}$ ( $D(A)
\cap D(B) \subset \Omega$ and $\overline{D(A) \cap D(B)} =
{\mathcal{H}}$). Thus $\xi$ is a closed densely defined sesquilinear
form. Assumptions $(H_{3})$ and $(H_{4})$ clearly imply that: there
exists a constant $const. = \max( c , c') > 0$ such that
$$ | \Im m \xi (u,u) | \; \leq \; const. \; \Re e \xi(u,u), \;
\; \forall u \in \Omega$$ Therefore $\xi$ is a sectorial
sesquilinear form. In summary, $\xi$ is a closed densely defined
sectorial form. According to Kato's first representation
theorem (see, e.g., [9, Theorem 2.1, pp. 322]): there exists a
unique m-sectorial operator associated with $\xi$ (m-sectorial extension of $\overline{A+B}$). Since
$\overline{A+B}$ is maximal and $\xi$ is
sectorial. Then $\overline{A+B}$ is m-sectorial, and it is the
m-sectorial operator associated with $\xi$. Since $D(A) = D(A^*)$
and $D(B) = D(B^*)$. It easily follows that $D(\overline{A+B})\subset
D((\overline{A+B})^*)$. Therefore
$D((\overline{A+B})^{\frac{1}{2}})\subset D((\overline{A+B})^{*
\frac{1}{2}})$. According to [10, Theorem 5.2, p. 238], we
have
\begin{equation}
D((\overline{A+B})^{\frac{1}{2}})\subset D(\xi) \subset
D((\overline{A+B})^{* \frac{1}{2}})
\end{equation}
\noindent In the same way, for the conjugate $\xi^*$ of $\xi$ we
have

\begin{equation}
D((\overline{A+B})^{* \frac{1}{2}})\subset D(\xi^*) \subset
D((\overline{A+B})^{\frac{1}{2}})
\end{equation}
Since $D(\xi) = D(\xi^*) = \Omega$ and from (2.1), (2.2). It follows
that $$ D((\overline{A+B})^{* \frac{1}{2}}) =  \Omega =
D((\overline{A+B})^{\frac{1}{2}})$$

We now consider our investigation related to the square
root problem of Kato for 
the "generalized" sum of operators defined above. We show
that the same conclusion still holds under appropriates assumptions.

We have

{\thm Let $A = A_{1} - i A_{2}$ and $B = B_{1} - i B_{2}$ be unbounded normal
operators on ${\mathcal H}$.
Assume that assumptions  $(H_{1})$, $(H_{3})$, $(H_{4})$, $(H_{5})$, and
$(H_{6})$ are satisfied.
Then there exists a unique m-sectorial operator
$A \oplus B$ satisfying the square root problem of Kato:
$$D((A \oplus B)^{\frac{1}{2}}) = \Omega(A) \cap \Omega(B) =
D((A \oplus B)^{ * \frac{1}{2}})$$ Also $A
\oplus B$ and $\overline{A+B}$ coincide if $\overline{A+B}$ is maximal.}

\vspace{0.4cm}

\noindent Proof. Consider the sesquilinear form $\xi = \phi +
\psi$ given by (1.4). Clearly $\xi$ is a
closed densely defined sequilinear form. Also since $\Re e
\xi(u,u) = \| (A_{1} \oplus B_{1})^{\frac{1}{2}} u \|^2 \; \;
\forall u \in \Omega $ and $\Im m \xi(u,u) = - \| ( A_{2} \oplus
B_{2})^{\frac{1}{2}} u\|^2 \; \; u \in \Omega$. It easily follows $\xi$ is
a sectorial sesquilinear form. Thus there exists
a $const = \max(c,c') > 0$ such that $| \Im m \xi(u,u) | \; \leq
\; const. \; \Re e \xi(u,u)$. According to Kato's first
representation theorem: there exists a unique m-sectorial
operator $A \oplus B$ associated with $\xi$ such that
$$\xi(u,v) = < \; (A \oplus B) u, v \; > \; \; \; u \in D(A
\oplus B), \; \; v \in \Omega$$
and in addition $D(A \oplus B) , \; \; D((A \oplus B)^*) \; \subset
D( \xi) = \Omega = \Omega_{\xi}$. Since $\Omega$ is
closed in $[ \Omega ,
{\mathcal{H}}]_{\frac{1}{2}}$. Then we conclude using a result
of Lions (see, e.g., [11, Theorem. 6.1, p. 238]), that
$A \oplus B$ satisfies the square root problem of Kato.
It is not hard to see that $A \oplus B = \overline{A + B}$ if
$\overline{A+B}$ is maximal. Indeed if $\overline{A+B}$ is maximal and since $A \oplus
B$ is an m-sectorial extension of it. Then they coincide everywhere.

\section{Applications}
In this section we give an application related to the algebraic
sum case. Consider the algebraic
sum given by a perturbation of the Schr\"{o}dinger operator
$S_{Z} = - Z \Delta + V$ with $Z$ a complex number and $V$
is a singular complex potential. Let $X \subset {\bf R}^{d}$ be an
open set and assume that our Hilbert space ${\mathcal{H}} =
L^2(X)$. Let $\Phi$ be the the sesquilinear form given by

\begin{equation}
\Phi_{Z}(u,v) = \int_{X} Z \nabla u \overline{\nabla v} dx, \; \; \;
\; \forall \; u,v \in D(\Phi) = {\bf H}_{0}^{1}(X)
\end{equation}
where $Z = \alpha - i \beta$ is a complex number satisfying the
following conditions: $ \alpha > 0, \; \; \beta > 0, \; \;
\mbox{and} \; \; \beta \leq \alpha$. The previous conditions on $Z$ clearly
that the sesquilinear form $\Phi$ is sectorial.

\noindent Let $V$ be a measurable complex function and let
$\Psi$ be the sesquilinear form given as
\begin{equation}
\Psi(u,v) = \int_{X} V u \overline{v} dx, \; \; \; \; \forall \; u,v
\in D(\Psi)
\end{equation}
with $D(\Psi) = \{ u \in L^2(X) \; : \; V | u |^2 \in L^1(X)
\}$. Throughout this section we assume that the potential $V \in
L_{loc}^1(X)$ and that there exists $\theta \in ( 0,
{\frac{\pi}{2}})$ such that
\begin{equation}
| \arg (V(x)) | \leq \theta, \; \; \; \mbox{almost everywhere}
\end{equation}
From (3.3) we have
\begin{equation}
| \Im m \Psi(u,u) | \leq \; \tan \theta \; \Re e \Psi(u,u), \; \; \;
\forall \; u \in D(\Psi)
\end{equation}
In other words, the sesquilinear form $\Psi$ is sectorial. Under
the previous assumptions $\Phi$ and $\Psi$ are
closed densely defined sectorial forms. The operators
associated with both $\Phi$
and $\Phi$ are
given as
$$D(A_{Z}) = \{ u \in {\bf H}_{0}^{1}(X) : Z \Delta u \in L^2(X)
\}, \; \; \; A_{Z} u = - Z \Delta u, \; \; \forall \; u \in
D(A_{Z})$$
$$D(B) = \{ u \in L^2(X) : V u \in L^2(X)\} , \; \; \; B u = V u, \;
\; \forall \; u \in D(B)$$ It is not hard to see that $A_{Z}$,
$B$ are unbounded normal operators in $L^2(X)$. One can
decompose them as $A_{Z} = A_{Z}^{1} - i A_{Z}^2$ where $A_{Z}^1
= - \alpha \Delta$ and $A_{Z}^2 = - \beta \Delta$ are
nonnegative selfadjoint operators, and $B = B_{V}^1 - i
B_{V}^2$  where $B_{V}^1$, $B_{V}^2$ are nonnegative
selfadjoint operators (that is justified by (3.3)). Now one assumes that
$X = {\bf R}^d$. Thus the operator associated with $\Xi = \Phi + \Psi$
(by the first representation theorem of Kato) is the closure of the
algebraic sum $S_{Z}$ (ie., $\overline{- Z \Delta + V}$) and it is
m-sectorial (see,  e.g., [3]). In fact such an operator
is defined as
$$D(\overline{- Z \Delta + V}) = \{ u \in {\bf H}^1({\bf R}^d) : V
|u|^2 \in L^1({\bf R}^d) \; \mbox{and} \; - Z \Delta u + V u \in
L^2({\bf R}^d) \}$$
$$\overline{- Z \Delta + V} u = -Z \Delta u + V
u, \; \forall u \in D(\overline{- Z \Delta + V})$$ Let us notice
that $D(A_{Z}) = H^2({\bf R}^d)$, $D(B) = \{ u \in L^2({\bf
R}^d) : V u \in L^2({\bf R}^d) \}$, and their intersection is dense
in $L^2({\bf R}^d)$. Therefore applying 
Theorem 2.1 to the operators $A_{Z}$ and $B$. It easily follows that
$$D((\overline{- Z \Delta + V})^{\frac{1}{2}}) = {\bf H}^1(R^d)
\cap D( B^{\frac{1}{2}}) = D((\overline{- Z \Delta + V})^{*
\frac{1}{2}})$$ For instance considering the case $d = 1$. Then we have
$$D((\overline{- Z \Delta + V})^{\frac{1}{2}}) = {\bf H}^1(R)
 = D((\overline{- Z \Delta + V})^{* \frac{1}{2}})$$

{\rem We may illustrate theorem 2.2 by considering both
$A_{Z}$ and $B$ defined above and by assuming that the potential $V \not = 0$
satisfies: (3.3) and the following
\begin{equation}
V \in L^1({\bf{R}^d}),  \; \; \; V \not \in L_{loc}^{2}({\bf{R}^d})
\end{equation}
In such a case, it is not hard to show that $D(A_{Z}) \cap D(B) =
\{ 0\}$ (see, e.g., [5, 8]). Therefore the algebraic sum $S_{z} = - Z \Delta +
V$ is not defined. Nevertheless, consider the sum $\Xi = \Phi +
\Psi$. Clearly $\Xi$ is a closed densely defined sectorial form
( $ C_{0}^{\infty} (\bf{R}^d) \subset D(\Phi) \cap D(\Psi)$).
According to the first representation theorem of Kato:
there exists a unique m-sectorial operator $( - Z \Delta
\oplus V)$ associated with $\Xi$. Thus theorem 2.2 can applied to
$A_{Z}$ and $B$. Therefore it easily follows the operator $( - Z \Delta \oplus V)$
satisfies the square root problem of Kato.}

\vspace{0.8cm}

\begin{center}
{\Large References}
\end{center}

\noindent [1] P. Ausher, S. Hofmann, A. McIntosh, and P.
Tchamitchian, {\it The Kato Square Root Problem for Higher Order
Elliptic Operators and Systems on ${\bf R}^n$}. J. Evol. Eq. {\bf
1} (2001), No. 4, pp. 361-385.

\vspace{0.2cm}

\noindent [2] A. Bivar-Weinholtz and M. Lapidus, {\it Product Formula
for Resolvents of Normal Operator and the Modified Feynman
Integral}. Proc. Amer. Math. Soc. Vol. {\bf 110}, No. 2 (1990).

\vspace{0.2cm}

\noindent [3] H. Br\'ezis and T. Kato, {\it Remarks on the
Schr\"{o}dinger Operator with Singular Complex Potentials}, J.
Math. Pures Appl. {\bf 58} (1979), 137-151.

\vspace{0.2cm}

\noindent [4] T. Diagana, {\it Sommes d'op\'erateurs et conjecture
de Kato-McIntosh}, C. R. Acad. Sci. Paris, t. {\bf 330}, S\'erie
I, p. 461-464 (2000).

\vspace{0.2cm}

\noindent [5] T. Diagana, {\it Schr\"{o}dinger operators with a
singular potential}. Int. J. Math-Math. Sci. Vol. {\bf 29}, No. 6
(2002), pp. 371-373.

\vspace{0.2cm}

\noindent [6] T. Diagana, {\it Quelques remarques sur
l'op\'erateur de Sch\"{o}dinger  avec un potentiel complexe
singulier particulier}. Bull. Belgian. Math. Soc. {\bf 9} (2002),
293 - 298.

\vspace{0.2cm}

\noindent [7] T. Diagana, {\it An application to Kato's square
root problem}. Int. J. Math-Math. Sci. Vol. {\bf 9}, No. 3 (2002),
pp. 179-181.

\vspace{0.2cm}

\noindent [8] T. Diagana, {\it  A Generalization related to
Schr\"{o}dinger operators with a singular potential}. Int. J.
Math-Math. Sci. Vol. {\bf 29}, No. 10 (2002), pp. 609-611.

\vspace{0.2cm}

\noindent [9] T. Kato, {\it Perturbation theory for linear
operators}, New york (1966).

\vspace{0.2cm}

\noindent [10] J. L. Lions, {\it Espace interm\'ediaires entre
espaces Hilbertiens et applications}. Bull. Math. Soc. Sci. Math.
Phys. R. P. Roumanie (N.S) {\bf 2} (50) (1958), 419-432.

\vspace{0.2cm}

\noindent [11] J. L. Lions, {\it Espaces d'interpolation et
domaines de puissances fractionnaires d'op\'erateurs}. J. Math.
Soc. Japan. {\bf 14} (2) (1962).

\vspace{0.2cm}

\noindent [12] A. Pazy, {\it Semigroups of linear operators and
application to partial differential equations}, Springer-Verlag,
New York (1983).

\vspace{0.2cm}

\noindent [13] W. Rudin, {\it Functional analysis}, Tata
McGraw-Hill, New Delhi (1974).

\end{document}